\documentclass[a4paper,12pt]{amsart}

\usepackage{amsmath,times,epsfig,amssymb,amsbsy,amscd,amsfonts,amstext,color,bm}

\usepackage[all]{xy}
\usepackage{bm}
\usepackage{comment}

\usepackage{tikz,epic,eso-pic}

\theoremstyle{plain}
\newtheorem{theorem}{Theorem}[section]
\newtheorem*{theorem*}{Theorem}
\newtheorem*{theoremA*}{Theorem A}
\newtheorem*{theoremB*}{Theorem B}
\newtheorem{corollary}[theorem]{Corollary}
\newtheorem{lemma}[theorem]{Lemma}

\theoremstyle{definition}
\newtheorem{definition}[theorem]{Definition}
\newtheorem{remark}[theorem]{Remark}
\newtheorem{example}[theorem]{Example}
\newtheorem{proposition}[theorem]{Proposition}

\newcommand{\Map}{\operatorname{Map}}
\newcommand{\Aut}{\operatorname{Aut}}
\newcommand{\Inn}{\operatorname{Inn}}

\newcommand{\id}{\operatorname{id}}
\newcommand{\ord}{\operatorname{ord}}

\newcommand{\bR}{\mathbb{R}}

\newcommand{\bZ}{\mathbb{Z}}

\newcommand{\cA}{\mathcal{A}}

\title{Homogeneous quandles with abelian inner automorphism groups}
\author{Takuya Saito, Sakumi Sugawara}
\date{\today}
\begin{document}
\maketitle

\begin{abstract}
In this paper, we give a characterization of homogeneous quandles with abelian inner automorphism groups. 
In particular, we show that such a quandle is expressed as an abelian extension of a trivial quandle.
Our construction is a generalization of the recent work by Furuki and Tamaru, which gives a construction of disconnected flat quandles.
\end{abstract}

\renewcommand{\thefootnote}{}
\footnotetext{
Takuya Saito: takuya.saito@unito.it\par
Sakumi Sugawara: sugawara.sakumi.f5@elms.hokudai.ac.jp
}
\renewcommand{\thefootnote}{\arabic{footnote}}

\section{Introduction}
A quandle is an algebraic structure introduced by Joyce \cite{joy} and Matveev \cite{mat}, where Matveev referred to them as distributive groupoids and they are widely used in the knot theory.
A \textit{quandle} is defined to be a pair of a set $Q$ and a \textit{point symmetry} $s_x:Q\to Q$ at each point $x\in Q$ and satisfies the following axioms;
\begin{itemize}
    \item [(Q1)] $s_x(x) = x$ for each $x \in Q$,
    \item [(Q2)] $s_x$ is bijective for each $x\in Q$, and
    \item [(Q3)] $s_{x} \circ s_{y} = s_{s_{x}(y)} \circ s_{x}$.
\end{itemize}
Each of these axioms corresponds to the Reidemeister move in knot diagrams.

By the way, Takasaki defined a kei in \cite{tak} as a generalization of symmetric transformations.
The axiom of keis is the axiom of quandles plus the condition $s_x^2=\id$ for all $x \in Q$. Hence, it is a specialization of quandles.
Symmetric spaces are also considered as quandles by ordinary point symmetry.
In this aspect, quandles are considered as a generalization of symmetric spaces.
Recently, some properties of quandles analogous to symmetric spaces have been studied, for example, homogeneity \cite{hig-kur}, doubly transitively \cite{tam}, antipodal subset \cite{kmtt}, and flatness \cite{ish-tam,fur-tam,sin}.

In particular, a flat symmetric space is a fundamental class of symmetric spaces, including tori.
Similar to a symmetric space, a quandle is \textit{flat} if the group generated by all compositions of two symmetries $s_x\circ s_y$ is abelian.
Ishihara and Tamaru classified finite flat connected quandles in \cite{ish-tam}. 
They showed that a finite quandle is flat and connected if and only if it is a discrete torus.
If quandles are not connected, they have not been classified.
Recently, Furuki and Tamaru found a method to construct a disconnected flat homogeneous quandle from a vertex-transitive graph in \cite{fur-tam}.
Let us recall their construction.
Let $\Gamma=(X,E)$ be a finite simple graph with the vertex set $X$ and the edge set $E$.
Define $d:X\times X\to \bZ_2:=\bZ/2\bZ$ by $d(x,y)=1$ if $x$ and $y$ are joined by an edge and $d(x,y)=0$ in other cases.
Let $Q$ denote $X\times \bZ_2$ and define a map $s_{-}(-):Q\to \Map(Q,Q)$ by $s_{(x,a)}(y,b)=(y,d(x,y)+b)$, where $\Map(X,X)$ denotes the set of maps from the set $X$ to itself.
Then, we have the following.
\begin{theorem*}[\cite{fur-tam}]
The quandle $Q_\Gamma := (Q, s)$ is disconnected and flat (moreover, its inner automorphism group is abelian). Furthermore, $Q_\Gamma$ is homogeneous if and only if the graph $\Gamma$ is vertex-transitive.
\end{theorem*}\noindent
Their construction produces a lot of disconnected flat quandles.
However, the complete classification of flat quandles is still a problem.

In this paper, we will focus on the homogeneous quandles with abelian inner automorphism groups.
If the inner automorphism group is abelian, then the quandle is flat.
Thus, our assumption is somewhat stronger than flatness.
We characterize homogeneous quandles with abelian inner automorphism groups.
Let $X$ be a set, $A$ be an abelian group and $d:X\times X \rightarrow A$ be a map satisfying $d(x,x) =0$ for all $x\in X$. 
Then, $Q=X \times A$ defines a quandle denoted by $X\times_{d} A$, with the binary operation $s_{(x,a)}(y,b) = (y,b+d(x,y))$ (see Section 3).
Our first result is the following.

\begin{theoremA*}[\textup{=Proposition \ref{prop:innner-comm}.+Theorem \ref{thm:main_transitive}.}]
When $(X,d)$ is homogeneous (in some sense), the quandle $X\times_{d} A$ is a homogeneous quandle with an abelian inner automorphism group.
\end{theoremA*}
Note that such a quandle is an abelian extension of a trivial quandle (see also \cite{cens,cks}).
Our construction is a generalization of Furuki--Tamaru's construction.
In fact, by specializing in $A=\mathbb{Z}_2$ in our construction, we have their theorem.
Our main result is the following.
\begin{theoremB*}[\textup{=Theorem 4.1.}]
Let $Q$ be a homogeneous quandle with an abelian inner automorphism group. Then, there exists a quandle $X \times_{d} A$ which is isomorphic to $Q$.
\end{theoremB*}
Using the above result, we will classify such quandles with small orders in Section \ref{sec:class}.
However, a pair $(X,d)$ is not necessarily homogeneous.
That is, the converse of Theorem A does not hold in general.
In Section \ref{sec:rema}, we construct an example that cannot be constructed from homogeneous $(X,d)$.

Other related studies on the commutativity of quandles are by Bardakov and Nasybullov in \cite{bar-nas}, Jedli\v{c}ka, Pilitowska, Stanovsk\'{y}, and Zamojska-Dzenio in \cite{jpsz} and by Lebed and Mortier in \cite{leb-mor}. However these do not assume homogeneity of quandles.
Bardakov and Nasybullov introduced the $(G,A)$-quandles, Jedli\v{c}ka et al. introduced the affine mesh to study medial quandles and Lebed--Mortier introduced the filtered-permutation to study quandles with abelian inner automorphism groups and abelian structure groups.
Our construction is also considered a specialization of these concepts.

The content of this paper is as follows.
Section \ref{sec:prel} contains the preliminaries on quandles and directed graphs. 
In Section \ref{sec:const}, we construct homogeneous quandles with abelian inner automorphism groups from directed graphs with edges weighted by an abelian group.
Section \ref{sec:char} shows that our construction provides all such quandles up to isomorphisms.
Section \ref{sec:class} gives a classification for small orders  and a formula for counting them of order $2p$ with prime $p$.
Section \ref{sec:rema} is a remark on Theorem A.
\medskip

\noindent
\textbf{Acknowledgment.}
The authors would like to thank Hiroshi Tamaru and Yuta Taniguchi for constructive discussion and helpful comments.
This work began with the intensive lecture on quandles by Professor Hiroshi Tamaru at Hokkaido University in November 2021.
This work is supported by JSPS KAKENHI Grant Numbers 22KJ0114 and 23KJ0031.

\section{Preliminary}\label{sec:prel}
\subsection{Basics on quandles}

We will write $\Map(X,X)$ as the set of maps from the set $X$ to itself.
Recall the definition of quandles.
\begin{definition}
Let $Q$ be a set and $s:Q \rightarrow \Map(Q,Q); x \mapsto s_{x} (-)$ be a map.
A pair $(Q,s)$ is called a {\it quandle} if 
\begin{itemize}
    \item [(Q1)] $s_x(x) = x$,
    \item [(Q2)] $s_x \in \Map(Q,Q)$ is bijective,
    \item [(Q3)] $s_{x} \circ s_{y} = s_{s_{x}(y)} \circ s_{x}$,
\end{itemize}
for each $x,y \in Q$.
\end{definition}

\begin{example}
The following pairs $(Q,s)$ are quandles;
\begin{itemize}    
    \item[(i)] Every set $X$ can have a quandle structure by defining the operation $s_{x}$ as the identity map $\mathrm{id}_{X}$ for any $x \in X$.
    This is called a \textit{trivial quandle}.
    \item[(ii)] Let $G$ be a group and define the operation $s:G \times G \rightarrow G$ as $s_{g} (h) = g^{-1} h g$. Then, $(G,s)$ is a quandle. This is called the conjugation quandle of a group.
    \item[(iii)] Let $\bR^{n}$ be a $n$-dimensional Euclidian space and define $s: \bR^{n} \times \bR^{n} \rightarrow \bR^{n}$ as $s_{x}(y) = 2x -y$. Then, the pair $(\bR^{n},s)$ is a quandle.
\end{itemize}
\end{example}

\begin{definition}
Let $(Q,s)$, $(Q',s')$ be quandles. A map $f : Q \rightarrow Q'$ is a quandle {\it homomorphism} if $f \circ s_{x} = s'_{f(x)} \circ f$ for each $x \in Q$. If $f$ is bijective, then $f$ is called a quandle {\it isomorphism}. 
\end{definition}

Note that a point symmetry $s_{x} : Q \rightarrow Q$ is a quandle isomorphism for any $x \in X$ by the axiom (Q3).

\begin{definition}
Let $(Q,s)$ be a quandle. 
An automorphism of $(Q,s)$ is a quandle isomorphism from $Q$ to itself. 
We denote the group of automorphisms by $\Aut(Q,s)$. This is called the {\it automorphism group} of a quandle $(Q,s)$.

Also, we denote by $\Inn(Q,s)$ the subgroup of $\Aut (Q,s)$ generated by the set $\{s_{x} \mid x \in Q\}$. This is called the {\it inner automorphism group}  of $(Q,s)$. 
\end{definition}

\begin{definition}
A quandle $(Q, s)$ is {\it homogeneous} (resp. {\it connected}) if $\Aut(Q, s)$ (resp. $\Inn(Q, s)$) acts transitively on $Q$.
\end{definition}

\subsection{Orbit decomposition by the inner automorphisms} 
Let $Q$ be a quandle and $Q = \bigsqcup_{x \in X} Qx$ be the orbit decomposition by the canonical $\Inn (Q)$-action. Here, $X = Q / \Inn(Q) $ denotes the orbit space, and each orbit $Qx$ is expressed as $\{g(x) \mid g \in \Inn(Q) \}$ for a $x \in Q$.

\begin{proposition}
Each orbit $Qx$ is a homogeneous quandle.
\end{proposition}

\proof
Since $Qx$ is $\Inn(Q)$-orbit, $s_{a}(b) \in Qx$ for each $a,b \in X$. Therefore, the binary operation is closed under $Qx$.
Thus, $Qx$ is a quandle. 
Moreover, since $\Inn(Q)$ acts on $Qx$ transitively, $Qx$ is a homogeneous quandle. 
\endproof

\begin{proposition} \label{prop:QxQy}
Let $Qx$ and $Qy$ be distinct orbits ($x,y \in X$). Suppose that there
exist $a \in Qx$ and $f \in \Aut(Q)$ such that $f(a) \in Qy$. Then, $f(Qx) \subset Qy$ and the restriction $f |_{Qx} : Qx \rightarrow Qy$ defines a quandle isomorphism.
\end{proposition}

\proof
Let $b \in Qx$. Since $\Inn(Q)$ acts on $Qx$ transitively, we can suppose that there exists $c \in Q$ such that $a = s_{c} (b)$. 
We have $f(a) = f \circ s_{c} (b) = s_{f(c)} \circ f(b) \in Qy$. Thus, it follows that $f(b) \in Qy$ since $Qy$ is an $\Inn(Q)$-orbit.
Therefore, $f(Qx) \subset Qy$ and $f|_{Qx}:Qx \rightarrow Qy$ is an injective quandle homomorphism. By the similar argument for $f^{-1}$, we have that $f|_{Qx}:Qx \rightarrow Qy$ is bijective. Thus, $f |_{Qx} : Qx \rightarrow Qy$ is a quandle isomorphism.
\endproof

\subsection{$A$-weighted digraph}
We recall some terminologies of directed graphs.

A directed graph (or just \textit{digraph}) is a pair $\Gamma=(X,D)$ of the \textit{vertex set} $X$, whose elements are called \textit{vertex}, and a subset $D$, called \textit{edge set}, of $X\times X$ whose elements are called \textit{directed edge}.
Two digraphs $(X,D)$ and $(X',D')$ are \textit{isomorphic} if there exists a bijection $f:X\to X'$ such that $(x,y)\in D$ if and only if $(f(x),f(y))\in D'$, and such a map is called a graph isomorphism or just an \textit{isomorphism}.
We write the set of all isomorphisms of a digraph $\Gamma$ by $\Aut (\Gamma)$.

Let $A$ be an abelian group.

\begin{definition}
    An $A$-\textit{weight} on a digraph $\Gamma=(X,D)$ is a map $d:D\to A$.
    We call a pair $(\Gamma,d)$ of a digraph and an $A$-weight an $A$-\textit{weighted graph}.
    An $A$-weight $d$ on $\Gamma$ is \textit{indecomposable} if $\{d(x,y)\in A\mid x\in X, (x,y)\in D\}$ generates $A$ as an abelian group for all $y\in X$.
\end{definition}

\begin{definition}
    Let $X$ be a set and let $d:X \times X \rightarrow A$ be a map satisfying $d(x,x) = 0$.
    Then, we can construct a digraph with the vertex set $X$ and the edge set $\{(x,y) \in X \times X \mid d(x,y)\neq 0\}$.
    This digraph is called the \textit{support graph} of $(X,d)$ and denoted by $\Gamma(X,d)$.    
\end{definition}

Note that a support digraph is loopless and has an $A$-weight induced by $d: X \times X \rightarrow A$.
From now on, an $A$-weighted graph defined in this way will be simply denoted by $(X,d)$, and $d$ is the $A$-weight on $(X,d)$.

\begin{definition}
    A graph isomorphism $f:X \rightarrow X$ is called \textit{$A$-isomorphism of $(X,d)$} if $d(f(x),f(y)) = d(x,y)$ for all $x,y \in  X$.
    We denote the set of all $A$-isomorphisms of $(X,d)$ by $\Aut (X,d)$. If $\Aut (X,d)$ acts on $X$ transitively, then we call $(X,d)$-\textit{homogeneous}.
\end{definition}

\section{The quandle $X\times_d A$}\label{sec:const}
We now define a quandle constructed by the direct product of a finite set and a group, which is the main object of this paper.
We define a set and a binary operation from an $A$-weighted graph $(X,d)$ as follows.
First, define the supporting set as $Q = X \times A$. 
The binary operation $s:Q \times Q \rightarrow Q$ is defined as $s_{(x,a)} (y,b) = (y, d(x,y) + b)$.
At this time, the following proposition holds.
\begin{proposition}
The pair $(Q,s)$ defined in this way is a quandle.
\end{proposition}

\proof
First, from the assumption $d(x,x) = 0$, $s_{(x,a)} (x,a) = (x,a)$ for any $(x,a) \in Q$. 
Second, the inverse of $s$ is defined as $s^{-1} : Q \times Q \rightarrow Q$ as $s^{-1}_{(x_a)} (y,b) = (y, -d(x,y)+b)$. Therefore, $s_{(x,a)}$ is bijective for any $(x,a) \in Q$.
Finally, we have
\[
s_{(x,a)} \circ s_{(y,b)} (z,c) = (z, d(x,z) + d(y,z)+b),
\]
 and 
 \[
 s_{s_{(x,a)}(y,b)} \circ s_{(x,a)} (z,c) = s_{(y, d(x,y)+b)} (z, d(x,z) + b) = (z, d(y,z) + d(x,z) + b).
 \]
 Therefore, the axiom (Q3) follows.
\endproof

We denoted by $(Q,s) = X\times_d A$ the quandle constructed in this way.

\begin{remark}
The quandle $X\times_d A$ is nothing but an abelian extension of a trivial quandle $X$ with an abelian group $A$ by a quandle $2$-cocycle $d: X \times X \rightarrow A$
(see also \cite{cens}, \cite{cks}).
\end{remark}

\begin{remark}
Since each $\Inn(X\times_d A)$-orbit is included by the set $\{x\}\times A$ for some $x\in X$,
the quandle $X\times_d A$ is disconnected unless $|X\times_d A|=1$ (see also Proposition \ref{prop:connected}).
\end{remark}

\begin{proposition}\label{prop:innner-comm}
The inner automorphism group of $X\times_d A$ is abelian.
\end{proposition}

\proof
It follows from $s_{(x,a)} \circ s_{(y,b)} (z,c) = s_{(y,b)} \circ s_{(x,a)}(z,c) = (z,c + d(x,z) + d(y,z))$ for each $(x,a) ,(y,b), (z,c) \in Q$.
\endproof

\begin{proposition}\label{prop:a-isom}
Let $f : (X,d) \rightarrow (X,d)$ be an $A$-isomorphism. The map $f' : Q \rightarrow Q$ defined as $f'(x,a) = (f(x),a)$ is a quandle isomorphism.
\end{proposition}

\proof
It is clear that $f'$ is bijective. 
It suffices to show that $f'$ is a quandle homomorphism.
For $(x,a), (y,b) \in Q$, we have
\[
f' \circ s_{(x,a)}(y,b) = f'(y,b+d(x,y)) = (f(y), b+d (x,y)),
\]
and 
\[
s_{f(x,a)} \circ f'(y,b) = s_{(f(x),a)} (f(y),b) = (f(y),b+d(f(x),f(y)) ) = (f(y), b+d (x,y)).
\]
Therefore, $f$ is a quandle isomorphism.
\endproof

\begin{proposition}
Let $\mu: X \rightarrow A$ be an arbitrary map. Then, the map $g_{\mu}: Q \rightarrow Q$ defined as $g_{\mu} (x,a) = (x,a + \mu(x))$ is a quandle isomorphism.
\end{proposition}

\proof
Since the map $g_{-\mu} : (x,a) \mapsto (a,a-\mu(x))$ gives the inverse, $g_{\mu}$ is bijective.
For each $(x,a), (y,b) \in Q$, we have 
$g_{\mu} \circ s_{(x,a)}(y,b) =s_{(x,a)}  \circ g_{\mu}(y,b) = (y, b+d(x,y) + \mu (x))$. Thus, $g_{\mu}$ is a quandle isomorphism.
\endproof

\begin{theorem}\label{thm:main_transitive}
If $(X,d)$ is $(X,d)$-homogeneous, then $X\times_d A$ is a homogeneous quandle. 
\end{theorem}

\proof
Let $(x,a), (y,b) \in X\times A$. By the assumption, there exists an $A$-isomorphism $f:X \rightarrow X$ such that $y = f(x)$. Let $\mu : X \rightarrow A$ be a map satisfying $\mu (x) = b-a$. Then,
$f' \circ g_{\mu} (x,a) = f'(x,b) = (y,b)$.
Therefore, $X\times_d A$ is a homogeneous quandle.
\endproof

By the above theorem, if $\Aut(X,d)$ acts on $X$ transitively, then $X\times_d A$ is particularly a homogeneous quandle with an abelian inner automorphism group.
Therefore, we can easily construct infinitely many homogeneous quandles whose inner automorphism groups are abelian.

Before considering the converse, let us discuss when $X\times_d A$ and $X\times_{d'} A$ are isomorphic for different $d$ and $d'$.
\begin{definition}
    Let $(X,d)$ be an $A$-weighted graph and let $\sigma:X\to\Aut(A)$ be a map assigning $y\in X$ to a group automorphism $\sigma_y:A\to A$.
    A \textit{flip} of $(X,d)$ by $\sigma$ is the $A$-weighted graph $(X,d^{\sigma})$ defined by:
    \begin{equation*}
        d^{\sigma}(x,y)=\sigma_y(d(x,y)).
    \end{equation*}
    Two $A$-weighted graphs $(X,d)$ and $(X',d')$ are \textit{flip equivalent} if there exists bijection $f:X\to X'$ and $\sigma:X\to\Aut(A)$ such that $d'(f(x),f(y))=d^\sigma(x,y)$.
    Such a pair $(f,\sigma)$ is called a \textit{weak isomorphism}.
\end{definition}
Clearly, flip equivalence is an equivalence relation.

\begin{proposition}\label{prop:flip}
    Let $f:X\to X$ be a bijection, let $d,d'$ be indecomposable $A$-weights, and let $g=(g_x:A\to A)_{x\in X}$ be a collection of bijections on $A$.
    A map $(f,g):X\times_d A\to X\times_{d'} A;(x,a)\mapsto (f(x), g_x(a))$ is a quandle isomorphism if and only if $(f,\sigma)$ is a weak isomorphism where $\sigma_x(a)=g_x(a)-g_x(0)$ for all $x\in X$ and $a\in A$.
    In particular, $f$ is an isomorphism between their support graphs.    
\end{proposition}

\proof
Firstly, assume $(f,g):X\times_d A\to X\times_{d'} A$ is a quandle isomorphism.
Let $s,s'$ denote operations on $X\times_{d} A$ and $X\times_{d'} A$, respectively.
Let $x,y\in X$ and $a,b\in A$.
Then
\begin{align*}
    &(f(y),g_y(d(x,y)+b))\\
    =&(f,g)(y,d(x,y)+b)\\
    =&(f,g)\circ s_{(x,a)}(y,b)\\
    =&s'_{(f,g)(x,a)} \circ (f,g)(y,b)\\
    =&s'_{(f(x),g_x(a))}(f(y),g_y(b))\\
    =&(f(y),d'(f(x),f(y))+g_y(b))
\end{align*}
holds.
Hence, we get $g_y(d(x,y)+b)=d'(f(x),f(y))+g_y(b)$, in other words, 
\begin{equation}\label{eq:iso}
    d'(f(x),f(y))=g_y(d(x,y)+b)-g_y(b).
\end{equation}
If $d(x,y)=0$, then $d'(f(x),f(y))=g_y(b)-g_y(b)=0$.
Thus, $f$ induces an isomorphism of the support graph.
The right side of the above equation \eqref{eq:iso} is a constant independent of $b$.
Substituting $0$ for $b$, $d'(f(x),f(y))=g_y(d(x,y))-g_y(0)$ holds.
Hence we get $g_y(d(x,y)+b)=g_y(d(x,y))+g_y(b)-g_y(0)$.
Define $\sigma_y:A\to A$ by $\sigma_y(b)=g_y(b)-g_y(0)$.
Then
\begin{align*}
    &\sigma_y(d(x,y)+b)=g_y(d(x,y)+b)-g_y(0)\\
    &=g_y(d(x,y))+g_y(b)-g_y(0)-g_y(0)=\sigma_y(d(x,y))+\sigma_y(b).
\end{align*}
By the assumption, $\{d(x,y)\mid x\in X\}$ spans $A$.
Thus $\sigma_y$ is an automorphism of $A$.
Therefore $d^\sigma(x,y)=g_y(d(x,y))-g_y(0)=d'(f(x),f(y))$ where $\sigma:x\mapsto \sigma_x$ and $(f,\sigma)$ is a weak isomorphism.
To prove the converse, follow the proof in reverse.
\endproof

\begin{corollary}\label{cor:q=2}
Suppose that $A=\bZ_2$. Then, quandles $X \times_d A$ and $X \times_{d'} A$ are isomorphic if and only if their support graphs $(X,d)$ and $(X,d')$ are isomorphic.
\end{corollary}

\proof
From Proposition \ref{prop:flip}, if quandles $X \times_d A$ and $X \times_{d'} A$ are isomorphic, then support graphs $(X,d)$ and $(X,d')$ are isomorphic.
Suppose that $(X,d)$ and $(X,d')$ are isomorphic. 
Then, there exists a bijective $f:X \rightarrow X$ such that $d(x,y) \neq 0$ if and only if $d(f(x),f(y))\neq 0$. 
Since we assume that $A=\bZ_2$, we have that $d(x,y) = d(f(x),f(y))$ for all $x,y \in X$. 
Therefore, from Proposition \ref{prop:a-isom}, it follows that quandles $X \times_d A$ and $X \times_{d'} A$ are isomorphic.
\endproof

\section{The characterization}\label{sec:char}
In this section, we characterize a homogeneous quandle with an abelian inner automorphism group. 
The following is the main result of the paper.

\begin{theorem}\label{thm:main}
Let $Q$ be a homogeneous quandle with an abelian inner automorphism group. Then, there exists a quandle $X\times_d A$ defined as the previous section which is isomorphic to $Q$.
\end{theorem}
We will give the precise statement in Theorem \ref{thm:main2}.
This theorem implies that the converse of Theorem \ref{thm:main_transitive} holds in some sense. However, the converse of Theorem \ref{thm:main_transitive} does not hold precisely (see Section 6).

Firstly, we see some properties of quandles with abelian inner automorphism groups.
Let $Q$ be a quandle and $Q = \bigsqcup_{x \in X} Qx$ be the $\Inn(Q)$-orbit decomposition.

\begin{proposition}\label{prop:equal}
Suppose that $\Inn(Q)$ is abelian and let $x \in X$. Then, $s_{a} = s_{b}$ for each $a,b \in Qx$. 
\end{proposition}

\proof
By the commutativity and the axiom (Q3), we have $s_{a} \circ s_{b} = 
s_{b} \circ s_{a} = s_{s_{b}(a)} \circ s_{b}$. Thus, $s_{a} = s_{s_{b}(a)}$ for each $a,b \in Q$. Hence, we have $s_{a} = s_{f(a)}$ for all $f \in \Inn(Q)$. Since $Qx$ is an $\Inn(Q)$-orbit, it follows that $s_{a} = s_{b}$.
\endproof

\begin{corollary}\label{cor:trivial}
Suppose that $\Inn(Q)$ is abelian. Then, each orbit $Qx$ is a trivial quandle.
\end{corollary}

\proof
Let $a,b \in Qx$. It follows from $s_{b}(a) = s_{a} (a) = a$.
\endproof

From here, throughout this section, we suppose that $Q$ is a homogeneous quandle and $\Inn(Q)$ is an abelian group.
Let $Q = \bigsqcup_{x \in X}Qx$ be an $\Inn(Q)$-orbit decomposition. 
Then, $Qx$ is isomorphic for all $x \in X$ from Proposition \ref{prop:QxQy}, 
and $Qx$ is a trivial quandle for all $x \in X$ for Corollary \ref{cor:trivial}.
Therefore, we have the following proposition.
\begin{proposition}
Let $a \in Q$ and $\Inn(Q)_{a} = \{g\in \Inn(Q) \mid g(a) = a\}$ be the
stabilizer group. Then, $\Inn(Q)_{a}$ is isomorphic for each $a \in Q$.
\end{proposition}

For each orbit, we fix an element $o_x \in Qx$. Let $A = \Inn (Q) / \Inn (Q)_{o_x}$. The quotient $A$ is an abelian group since $\Inn(Q)$ is abelian. By the above proposition, the definition of $A$ is independent of the choice of $x$ and $o_x$ up to isomorphism. 
Define a map $\varphi_{x} : A \rightarrow Qx $ as $\varphi([g]) = g(o_{x})$ for $g \in \Inn(Q)$. The symbol $[g]$ stands for the equivalent class of $g$ in $A$. By the orbit stabilizer theorem, the map $\varphi_{x}$ is well-defined and bijective. 

Note that $\varphi_x(0) = o_x$ and $\varphi_x ([g]+[h]) = g(h(o_x))$ for $g,h \in \Inn(Q)$, where $0$ is the identity element of $A$ and $+$ is the group operation in $A$. 

Since $s_{a} = s_{b}$ for each $a,b \in Qx$ by Proposition \ref{prop:equal}, we denote it by $s_{x}$.

\begin{definition}
For $x,y \in X$, we define $d: X \times X \rightarrow A$ as 
\[
d(x,y) = \varphi_{y}^{-1} (s_{x} (o_y)).
\]
\end{definition}
Since $Qx$ is a trivial quandle, $d(x,x) = 0$ for every $x \in X$.
We can define a quandle $Q' = X\times_d A$ by using the above settings, that is, the supporting set $Q' = X \times A$ and the binary operation $s'_{(x,a)} (y,b) = (y, d(x,y)+b)$.

With the preceding setting, we have the following theorem, which also establishes our main result.
\begin{theorem}\label{thm:main2}
Define a map $f: Q \rightarrow Q'$ as $f(\alpha) = (x,\varphi_{x}^{-1} (\alpha))$. Here $\alpha \in Qx$ and $\varphi_x$ is the map defined before. Then, $f$ is a quandle isomorphism. \end{theorem}

\proof
Since the map $(x,a) \mapsto \varphi_{x} (a) \in Qx \subset Q$ defines the inverse, $f$ is a bijection.
We will show that $f$ is a quandle homomorphism.
Let $\alpha, \beta \in Q$ and suppose that $\alpha \in Qx$ and $\beta \in Qy$.
Similarly, as before, we denote $s_{\alpha}$ by $s_{x}$. Note that there exists $b \in A$ such that $\beta= b (o_y)$. Then, we have
\begin{eqnarray*}
\varphi_{y}^{-1} (s_{x}(\beta)) &=&\varphi_{y}^{-1} (s_{x}(b(o_y))) \\
&=& \varphi_{y}^{-1} (b(s_{x}(o_y))) \\
&=& \varphi_{y}^{-1} (s_{x}(o_y)) + b \\
&=& d(x,y) + b.
\end{eqnarray*}
Thus, $f \circ s_{x} (\beta) = (y, \varphi_{y}^{-1}(s_{x}(\beta))) = (y, d(x,y)+b)$. 

On the other hand, we have 
\begin{eqnarray*}
s'_{f(\alpha)} \circ f (\beta) &=& s'_{(x,\varphi_{x}^{-1}(\alpha))} (y, \varphi_{y}^{-1}(\beta)) \\
&=& s'_{(x,\varphi_{x}^{-1}(\alpha))} (y, b) \\
&=& (y, d(x,y)+b).
\end{eqnarray*}
Therefore, $f$ defines a quandle isomorphism.
\endproof

Therefore, every homogeneous quandle with an abelian inner automorphism group is constructed by the pair of $X=Q / \Inn(Q)$ and the map $d: X \times X \rightarrow A$ with $d(x,x)=0$ for all $x \in X$.
Thus, we obtain an $A$-weighted graph $(X,d)$ and call it a \textit{digraph presentation} of $Q$.

If a digraph presentation $(X,d)$ of $Q$ is not indecomposable, then the number of $\Inn(Q)$-orbits of $\{x\}\times A$ should be more than one.
Hence, we obtain the following proposition.
\begin{proposition}
    For any $Q$, every digraph presentation of $Q$ is an indecomposable $A$-weighted graph.
\end{proposition}

\begin{remark}
In \cite{leb-mor}, quandles with abelian inner automorphism groups are characterized in terms of a filtered-permutation quandle.
For each $k\in\{1,\cdots, r\}$, let $M^{(k)}=(m^{(k)}_{i,j})_{1\leq i<j\leq r_k}$ be a collection of $r(r-1)/2$ integers with the conditions $1\leq m_{i,i}$ and $0\leq m_{i,j}<m_{j,j}(i\leq j)$.
Define $G(M^{(k)})$ as an additive group 
$\langle x(1,k),\dots,x(r-1,k)\mid \sum_{j=1}^{i}m^{(k)}_{i,j}x(j,k)=0\rangle$ of order $m^{(k)}_{1,1}m^{(k)}_{2,2}\dots m^{(k)}_{r-1,r-1}$.
The filtered-permutation quandle $Q(M)=\allowbreak (Q(M^{(1)},\dots,M^{(r)}),s)$ is defined on the disjoint union $Q(M^{(1)},\dots,M^{(r)})=\bigsqcup_{k=1}^r G(M^{(k)})$ with a point symmetry $s$ defined by $s_{b}(a)=a+x(j-i,i)\in M^{(i)}$ for $a\in M^{(i)}$ and $b\in M^{(j)}$.
By definition, each $G(M^{(k)})$ is an $\Inn(Q(M))$-orbit.
Hence, if $Q(M)$ is homogeneous, then $G(M)\cong G(M^{(k)})$ is independent of $k$ up to isomorphisms and $Q(M)$ is written by $\{1,\dots,r\}\times_d G(M)$ where $d(i,j)=x(j-i,i)$.
This agrees with our presentation.
\end{remark}

\begin{remark}
In \cite{jpsz}, medial quandles are characterized in terms of the sum of an indecomposable affine mesh.
A quandle $Q$ is medial if $s_{s_{a}(b)}\circ s_{c}=s_{s_{a}(c)}\circ s_{b}$ for all $a,b,c\in Q$.
If $\Inn(Q)$ is abelian, then we get
\begin{equation*}
    s_{s_{a}(b)}\circ s_{c}=s_{b}\circ s_{c}=s_{c}\circ s_{b}=s_{s_{a}(c)}\circ s_{b}
\end{equation*}
since $s_{s_{a}(b)}=s_b$ for all $a,b\in Q$.
Hence $Q$ is a medial quandle.
We can characterize a homogeneous quandle with an abelian inner automorphism group in terms of an affine mesh.

An affine mesh over a non-empty set $X$ is the triple $$((A_x)_{x\in X},(\phi_{x,y})_{x,y\in X},(d(x,y))_{x,y\in X})$$ with the following conditions;
\begin{enumerate}
    \item $A_x$ is an abelian group for each $x\in X$;
    \item $\phi_{x,y}:A_x\to A_y$ is a group homomorphism;
    \item $d(x,y)\in A_y$ and $d(x,x)=0$;
    \item $1-\phi_{x,x}\in\Aut(A_x)$;
    \item $\phi_{y,z}\phi_{x,y}=\phi_{y',z}\phi_{x,y'}$, and;
    \item $\phi_{y,z}(d(x,y))=\phi_{z,z}(d(x,z)-d(y,z))$.
\end{enumerate}
The sum of an affine mesh $\cA=((A_x)_{x\in X},(\phi_{x,y})_{x,y\in X},(d(x,y))_{x,y\in X})$ is a quandle on $\bigsqcup_{x\in X}A_x$ with point symmetry $s$ defined as 
\begin{equation*}
    s_a(b)=b+d(x,y)+\phi_{x,y}(a)-\phi_{y,y}(b)
\end{equation*}
for $a\in A_x,b\in A_y$.
The affine mesh $\cA=((A_x)_{x\in X},(\phi_{x,y})_{x,y\in X},(d(x,y))_{x,y\in X})$ is indecomposable if and only if $d(x,y)$ and $\phi_{x,y}(a)$ with $x\in X,a\in A_x$ span $A_y$ for every $y\in X$.
If the affine mesh is indecomposable, then $A_x$ is $\Inn(\cA)$-orbit for any $x\in X$ (see \cite[Lemma 3.10]{jpsz}).

If the sum of indecomposable affine mesh $\cA$ is homogeneous, all $A_x$ are isomorphic to each other for all $x\in X$ because each $A_x$ is an $\Inn(\cA)$-orbit.
In addition, suppose that $\cA$ has an abelian inner automorphism group, then $\phi_{x,y}(a)-\phi_{y,y}(b)$ depends on only $x,y\in X$.
Since both $\phi_{x,y}$ and $\phi_{y,y}$ are homomorphisms, we have $\phi_{x,y}(a)=\phi_{y,y}(b)=0$.
Therefore, $\cA$ is isomorphic to $X\times_{d'} A_{x_0}$ where $d'(x,y)$ is defined as the image of isomorphisms $A_y\to A_{x_0}$ of $d(x,y)$ for $x,y\in X$.
This can be considered as our presentation.
\end{remark}

\section{Classification for small orders}\label{sec:class}
In this section, we classify homogeneous quandles with abelian inner automorphism groups for small orders.
Throughout this section, we suppose that $Q$ is homogeneous, $\Inn(Q)$ is abelian, and $|Q| = n$ is finite.
Also, we identify $X$ with $Q/\Inn(Q)$.
The cardinality $|Q|$ is decomposed as $n=|X| \cdot |A|$.

In the following Table \ref{table}, the top row is the order $n$, and the bottom row is the number of isomorphism classes of homogeneous quandles of order $n$ with abelian inner automorphism groups.
\begin{table}[h]
\begin{equation*}
	\begin{array}{l||c|c|c|c|c|c|c|c|c|c|c|c|c|c|c|l}
	\mbox{order }&1&2&3&4&5&6&7&8&9&10&11&12&13&14&15&\cdots \\ \hline
	\sharp&1&1&1&2&1&4&1&7&4&7&1&36&1&15&18&\cdots
	\end{array}
\end{equation*}
\caption{The number of isomorphism classes}\label{table}
\end{table}

\subsection{General cases and tools}
\begin{proposition}\label{prop:connected}
Suppose that $Q$ is connected. Then $Q$ is a trivial quandle of a singleton.
\end{proposition}

\proof
Since $Q$ is connected, $Q$ is the $\Inn(Q)$-orbit itself. From Corollary \ref{cor:trivial}, $Q$ is a trivial quandle. 
Since the $\Inn(Q)$-orbit of a trivial quandle is a singleton, $Q$ is a singleton.
\endproof

\begin{proposition}\label{prop:primesize}
Suppose that $|Q| = p$ is a prime number. Then, $Q$ is a trivial quandle.
\end{proposition}

\proof
Since the cardinality of each orbit $|Qx|$ is constant for all $x \in X$, $|Q| = |X| \cdot |Qx|$. Thus, $|X|=1, |Qx|=p$ or $|X|=p, |Qx|=1$. If $|X| = 1$, then $Q$ is an $\Inn(Q)$-orbit. This does not happen by Proposition \ref{prop:connected}. Therefore, $|X|=p$ and $|Qx|=1$. Since every $\Inn(Q)$-orbit is a singleton, we have that $Q$ is a trivial quandle.
\endproof

\begin{proposition}\label{prop:2orbits}
Suppose that $\Inn(Q)$-orbit consists of two elements $X= \{x,y\}$. Then, $Q$ is isomorphic to $Q(\{x,y\}, \bZ_{\frac{n}{2}}, d)$, 
where $d(x,x) = d(y,y) = 0$ and $d(x,y) = d(y,x) = 1$ (see Figure \ref{fig:2orbit}).
\end{proposition}

\proof
The orbit $Qx$ is expressed as $\{g(o_x) \mid g \in \Inn(Q) \}$. 
Since there are only two $\Inn(Q)$-orbits and $s_{a}(o_{x}) = o_{x}$ for $a \in Qx$, we have 
\[
Qx = \{s_{b}^{k} (o_x) \mid b \in Qy, k \in \bZ \}.
\]
Each $s_{b}$ defines the same map independent of the choice of $b \in Qy$, which is denoted by $s_{y}$. Therefore, $\varphi^{-1} (Qx) $ is a cyclic group generated by $\varphi^{-1}_{x} (s_{y} (o_x))$. Thus $A \cong \bZ_{\frac{n}{2}}$. By identifying $1 \in \bZ_{\frac{n}{2}}$ with $\varphi^{-1}_{x} (s_{y} (o_x))$. We have $d(x,y) = d(y,x) = 1$.
\endproof

\begin{figure}[h]
    \centering
    \begin{tikzpicture}[auto]
	\node[draw,circle] (a) at (-1,0){$x$};
	\node[draw,circle] (b) at (1,0){$y$};

	\draw[->,thick,>=stealth] (a) to[bend right=30] node[auto=right]{$+1$} (b);
	\draw[->,thick,>=stealth] (b) to[bend right=30] node[auto=right]{$+1$} (a);
	\end{tikzpicture}
    \caption{$|X|=2$}
    \label{fig:2orbit}
\end{figure}

From these propositions, the isomorphism class of a quandle $Q$ is uniquely determined when $|X| = 1,2,n$. Next, we consider the case $|X| = p$ is a prime.

The following is a general proposition from group theory.
\begin{proposition}\label{prop:cyclic}
Let $X$ be a finite set whose cardinality is divided by a prime $p$, and suppose that a finite group $G$ acts on $X$ transitively. Then, $G$ contains the cyclic group of order $p$ as a subgroup.
\end{proposition}

\proof
Fix an element $x \in X$ and let $G_{x}$ be the stabilizer subgroup. Since $G$ acts transitively, $|G| = |X| \cdot |G_x|$. Since $p \mid |X|$, we get $p \mid |G|$. Thus, there exists an element of order $p$ by Cauchy theorem. The subgroup generated by it gives the desired one.
\endproof

Let $G_Q = \{\varphi \in \Aut(Q) \mid \varphi(x,a) \in \{x\} \times A \mbox{ for each }x \in X\mbox{ and for each }a\in A \}.$ Then, $G_Q$ is a normal subgroup of $\Aut(Q)$. Let $H_Q = \Aut(Q) / G_Q$ be the quotient group. Note that $H_Q$ acts on $X$ transitively.

\begin{proposition}\label{prop:primeaction}
Suppose that $|X|=p$ for $Q = X \times A$. Then, there exists a subgroup of $\Aut(Q)$ which is isomorphic to $\bZ_p:=\bZ/p\bZ$ and acts on $X$ transitively.
\end{proposition}

\proof
By Proposition \ref{prop:cyclic}, $H_Q$ contains $\bZ_p$ as a subgroup and let $\psi$ be a generator. Then, $\psi$ gives a non-trivial action on $X$, so there exist distinct $x,y \in X$ such that $y = \psi(x)$. 
The order of $\psi$ in $H_Q$ and $|X|$ are both $p$, and we have that $X = \{\psi^{k}(x) \mid 0 \leq k \leq p-1\}$. 
We can choose an element $\overline{\psi} \in \Aut(Q)$ as the lift of $\psi \in H$ defined as $\overline{\psi} (x ,a) = (\psi(x),a)$. Thus, the group generated $\overline{\psi}$ gives the desired one.
\endproof

\begin{remark}
By Proposition \ref{prop:primeaction}, we can identify $X$ with $\bZ_p$ as $\bZ_p$-set with the action as in the proof. $k \in \bZ_p$ is identified with $\psi^{k} (0)$. 
\end{remark}

\begin{proposition}\label{prop:primesize2}
Suppose that $|X|=p$ is a prime and fix $y_0 \in X$. The quandle $Q=X\times_d A$ is determined by the list $(d(x,y_0))_{x\in X}$.
\end{proposition}
\proof
From the above proposition, $X$ is expressed as $\{\psi^{k}(y_0) \mid 0\leq k \leq p-1\}$. 
Since the map $\overline{\psi}$ is a quandle homomorphism, we have $d(x,y)=d(\psi(x), \psi(y))$.
For each $y \in X$, there exists $k_0$ such that $y = \psi^{k_0} (y_0)$. 
Thus we have $d(x,y) = d(\psi^{-k_0}(x), y_0)$ for $x,y \in X$.
The list $(d(x,y_0))_{x\in X}$ determines $(d(x,y))_{x,y \in X}$.
\endproof

\begin{proposition}\label{prop:cyclic}
If there is $f\in H_Q$ whose order is equal to $|X|$, then $f$ induces a self $A$-isomorphism on the digraph presentation $(X,d)$.
\end{proposition}
\proof
Since the order of $f$ is equal to $|X|$, $f$ acts transitively on $X$, and $X$ can be expressed as $\{f^{k} (x) \mid 0\leq k \leq |X|-1\}$.
Similarly, as above, we can take the lift $\overline{f} \in \Aut(Q)$ of $f$ as $\overline{f}(x,a) = (f(x),a)$. 
Since $\overline{f}$ is a quandle homomorphism, we have $d(x,y)=d(f(x),f(y))$. 
Thus, we have the proposition.
\endproof

Finally, before listing up the isomorphism classes of homogeneous quandles $X\times_d A$, let us see how to determine isomorphisms between support graphs.

\begin{definition}
    Let $m$ be a positive integer and let $S$ be a subset of $\bZ_m$.
    A circulant graph $\mathrm{Cay}(\bZ_m,S)$ is a digraph with the vertex set $\bZ_m$ and the edge set $\{(i,j)\in\bZ_m\times\bZ_m \mid j-i\in S\}$.
\end{definition}
\begin{proposition}
If there is $f\in H_Q$ whose order is equal to $|X|$, then there exists a digraph presentation $(X,d)$ such that $\Gamma(X,d)$ is isomorphic to a circulant graph $\mathrm{Cay}(\bZ_{|X|},S)$ for some $S\subset\bZ_{|X|}$.
\end{proposition}
\proof
By Proposition \ref{prop:flip}, $f$ is a graph isomorphism.
Fix $x_0\in X$ and let $S$ be $\{i\in \bZ_{|X|}\mid  d(f^i(x_0),x_0)\neq 0\}$.
Then we get
\begin{align*}
    &j-i\in S\\
    \Leftrightarrow&d(f^{j-i}(x_0),x_0)\neq0\\
    \Leftrightarrow&d(f^{j}(x_0),f^{i}(x_0))\neq0\\
    \Leftrightarrow&(f^{j}(x_0),f^{i}(x_0))\mbox{ is a edge of }\Gamma(X,d).
\end{align*}
Therefore, $\Gamma (X,d)$ is isomorphic to $(\bZ_{|X|}, S)$.
\endproof
Muzychuk characterized the isomorphism class of circulant graphs in \cite{muz}.
In particular, for a prime $p$, it holds that two circulant graphs $\mathrm{Cay}(\bZ_{p},S)$ and $\mathrm{Cay}(\bZ_{p},S')$ are isomorphic if and only if there exists $\sigma\in\Aut(\bZ_p)$ such that $S'=\{\sigma.s\mid s\in S\}$.
Hence, when $|X|=p$ is a prime number, if we identify $X$ with $\bZ_p$ via the bijection $\bZ_p\to Q/\Inn(Q);i\mapsto f^i(x_0)$ for fixed $f\in H_Q\setminus\{\id\}$, then an isomorphism $\bZ_p\times_d A\to \bZ_p\times_d A$ induces an automorphism of $\bZ_p$.

\subsection{Each case}
Let us list the equivalent class of homogeneous quandles of order {$10$} or less with abelian inner automorphism groups.

By Proposition \ref{prop:primesize}, we know such quandles of prime order are trivial.
The second simplest case is $|X\times_d A|=pq$, where both $p$ and $q$ are prime.
In this case, the number of $\Inn(Q)$-orbits can be $p,q$, or $pq$.
If it is $pq$, then $Q$ is trivial.
Suppose $|X|=p$.
By Proposition \ref{prop:primeaction}, there exists $f\in H_Q$ of the order $p$.
From Proposition \ref{prop:primesize2}, $Q$ is determined by lists of weights
\begin{equation*}
    (d(x_1,x_0),\dots,d(x_{p-1},x_0)) (\neq (0,\dots,0))
\end{equation*}
where $x_0\in X$ and $x_k=f^k(x_0)$ for $0\leq k \leq p-1$.
We can suppose that $d(x_1,x_0)=1$ by $\bZ_p$ and $\Aut(\bZ_q)$ action.
Thus, we choose representatives of lists of weights of the following form
\begin{equation*}
    (d(x_1,x_0),\dots,d(x_{p-1},x_0))=(1,d_2,\dots,d_{p-1}).
\end{equation*}
where $d_i\in \bZ_q$ for each $i\in\{2,\dots,p-1\}$.

Let us consider the cases $n=2 \times \mbox{(prime)}$ ($|Q|=n=4,6,10$).
When $|X|=2$, there is only one isomorphism class from Proposition \ref{prop:primesize}.
When $|A|=2$, it suffices to compute the isomorphism class of digraphs from Corollary \ref{cor:q=2}.
Thus, if $n=4$, there are only two isomorphism classes; each class corresponds to the case $|X|=2,4$.

Next, we will compute the case when $n=6,10$.

\begin{example}[$n=6$]
When $|X|=3$, possible representatives of lists of weights are $(1,0)$ and $(1,1)$.
They give different support graphs.
The support graph of $(X,d)$ is an invariant.
Therefore, we get the list, and the digraph presentations of non-trivial quandles are in Figure \ref{fig:order6}.
Here, each $(x_{i},x_{j})$ entry stands for $d(x_{i},x_{j})$.
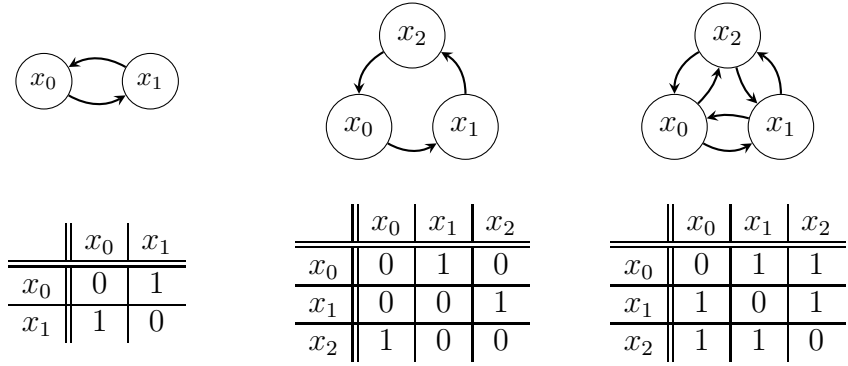
\begin{figure}[h]
    \centering
    \begin{tabular}{ccc}
    \begin{minipage}{0.3\linewidth}
    \centering
    \begin{tikzpicture}[auto,scale=0.7]\footnotesize
	\node[draw,circle] (a) at (-1,0){$x_0$};
	\node[draw,circle] (b) at (1,0){$x_1$};
	\draw[->,thick,>=stealth] (a) to[bend right=30] (b);
	\draw[->,thick,>=stealth] (b) to[bend right=30] (a);
	\end{tikzpicture}
    \end{minipage}&
    \begin{minipage}{0.3\linewidth}
    \centering
    \begin{tikzpicture}[auto,scale=0.7]
	\node[draw,circle] (a) at (-1,0){$x_0$};
	\node[draw,circle] (b) at (1,0){$x_1$};
	\node[draw,circle] (c) at (0,1.73){$x_2$};

	\draw[->,thick,>=stealth] (a) to[bend right=30] (b);
	\draw[->,thick,>=stealth] (b) to[bend right=30] (c);
	\draw[->,thick,>=stealth] (c) to[bend right=30] (a);
	\end{tikzpicture}
    \end{minipage}&
    \begin{minipage}{0.3\linewidth}
    \centering
    \begin{tikzpicture}[auto,scale=0.7]
	\node[draw,circle] (a) at (-1,0){$x_0$};
	\node[draw,circle] (b) at (1,0){$x_1$};
	\node[draw,circle] (c) at (0,1.73){$x_2$};

	\draw[->,thick,>=stealth] (a) to[bend right=30] (b);
	\draw[->,thick,>=stealth] (b) to[bend right=30] (c);
	\draw[->,thick,>=stealth] (c) to[bend right=30] (a);

	\draw[->,thick,>=stealth] (a) to[bend right=15] (c);
	\draw[->,thick,>=stealth] (b) to[bend right=15] (a);
	\draw[->,thick,>=stealth] (c) to[bend right=15] (b);
	\end{tikzpicture}
    \end{minipage}\\
    \vspace{1mm}\\
    \mbox{$\begin{array}{r||c|c}
	&x_0&x_1 \\ \hline\hline
	x_0&0&1 \\ \hline
	x_1&1&0
	\end{array}$}&
    \mbox{$\begin{array}{r||c|c|c}
	&x_0&x_1&x_2 \\ \hline\hline
	x_0&0&1&0 \\ \hline
	x_1&0&0&1 \\ \hline
	x_2&1&0&0
	\end{array}$}&
    \mbox{$\begin{array}{r||c|c|c}
	&x_0&x_1&x_2 \\ \hline\hline
	x_0&0&1&1 \\ \hline
	x_1&1&0&1 \\ \hline
	x_2&1&1&0
	\end{array}$}
    \end{tabular}
    \caption{The case $n=6$}
    \label{fig:order6}
\end{figure}
\end{example}

\begin{example}[$n=10$]\label{exam:n=10}
When $|X|=5$, the possible representatives of weight lists are $(1,0,0,0),(1,0,0,1),(1,0,1,0), (1,0,1,1),\allowbreak(1,1,0,0),(1,1,0,1),(1,1,1,0)$, and $(1,1,1,1)$.

These weights are classified into five types as follows
\begin{itemize}
    \item[(1)] $(1,0,0,0)$;
    \item[(2a)] $(1,0,1,0),(1,1,0,0)$;
    \item[(2b)] $(1,0,0,1)$;
    \item[(3)] $(1,0,1,1),(1,1,0,1),(1,1,1,0)$;
    \item[(4)] $(1,1,1,1)$.
\end{itemize}
Furthermore, they give different support graphs, as shown in Figure \ref{fig:order10}.
\begin{figure}[h]
    \centering
    \begin{tabular}{ccccc}
    \begin{minipage}{0.15\linewidth}
        \centering
        \begin{tikzpicture}[auto,scale=0.7,rotate=0]
			\node[draw,circle,inner sep=2pt] (a) at (0,1){};
			\node[draw,circle,inner sep=2pt] (b) at (sin{-72}, cos{-72}){};
			\node[draw,circle,inner sep=2pt] (c) at (sin{-144}, cos{-144}){};
			\node[draw,circle,inner sep=2pt] (d) at (sin{-216}, cos{-216}){};
			\node[draw,circle,inner sep=2pt] (e) at (sin{-288}, cos{-288}){};
			\draw[->,thick,>=stealth] (a) to[bend right=36] (b);
			\draw[->,thick,>=stealth] (b) to[bend right=36] (c);
			\draw[->,thick,>=stealth] (c) to[bend right=36] (d);
			\draw[->,thick,>=stealth] (d) to[bend right=36] (e);
			\draw[->,thick,>=stealth] (e) to[bend right=36] (a);
			\end{tikzpicture}
            \mbox{(1)}
    \end{minipage}&
    \begin{minipage}{0.15\linewidth}
        \centering
        \begin{tikzpicture}[auto,scale=0.7,rotate=0]
			\node[draw,circle,inner sep=2pt] (a) at (0,1){};
			\node[draw,circle,inner sep=2pt] (b) at (sin{-72}, cos{-72}){};
			\node[draw,circle,inner sep=2pt] (c) at (sin{-144}, cos{-144}){};
			\node[draw,circle,inner sep=2pt] (d) at (sin{-216}, cos{-216}){};
			\node[draw,circle,inner sep=2pt] (e) at (sin{-288}, cos{-288}){};
			\draw[->,thick,>=stealth] (a) to[bend right=36] (b);
			\draw[->,thick,>=stealth] (b) to[bend right=36] (c);
			\draw[->,thick,>=stealth] (c) to[bend right=36] (d);
			\draw[->,thick,>=stealth] (d) to[bend right=36] (e);
			\draw[->,thick,>=stealth] (e) to[bend right=36] (a);
			\draw[->,thick,>=stealth] (a) to[bend right=15] (c);
			\draw[->,thick,>=stealth] (b) to[bend right=15] (d);
			\draw[->,thick,>=stealth] (c) to[bend right=15] (e);
			\draw[->,thick,>=stealth] (d) to[bend right=15] (a);
			\draw[->,thick,>=stealth] (e) to[bend right=15] (b);
			\end{tikzpicture}
            \mbox{(2a)}
    \end{minipage}&
    \begin{minipage}{0.15\linewidth}
        \centering
        \begin{tikzpicture}[auto,scale=0.7,rotate=0]
			\node[draw,circle,inner sep=2pt] (a) at (0,1){};
			\node[draw,circle,inner sep=2pt] (b) at (sin{-72}, cos{-72}){};
			\node[draw,circle,inner sep=2pt] (c) at (sin{-144}, cos{-144}){};
			\node[draw,circle,inner sep=2pt] (d) at (sin{-216}, cos{-216}){};
			\node[draw,circle,inner sep=2pt] (e) at (sin{-288}, cos{-288}){};
			\draw[<->,thick,>=stealth] (a) to[bend right=36] (b);
			\draw[<->,thick,>=stealth] (b) to[bend right=36] (c);
			\draw[<->,thick,>=stealth] (c) to[bend right=36] (d);
			\draw[<->,thick,>=stealth] (d) to[bend right=36] (e);
			\draw[<->,thick,>=stealth] (e) to[bend right=36] (a);
			\end{tikzpicture}
            \mbox{(2b)}
    \end{minipage}&
    \begin{minipage}{0.15\linewidth}
        \centering
        \begin{tikzpicture}[auto,scale=0.7,rotate=0]
			\node[draw,circle,inner sep=2pt] (a) at (0,1){};
			\node[draw,circle,inner sep=2pt] (b) at (sin{-72}, cos{-72}){};
			\node[draw,circle,inner sep=2pt] (c) at (sin{-144}, cos{-144}){};
			\node[draw,circle,inner sep=2pt] (d) at (sin{-216}, cos{-216}){};
			\node[draw,circle,inner sep=2pt] (e) at (sin{-288}, cos{-288}){};
			\draw[<->,thick,>=stealth] (a) to[bend right=36] (b);
			\draw[<->,thick,>=stealth] (b) to[bend right=36] (c);
			\draw[<->,thick,>=stealth] (c) to[bend right=36] (d);
			\draw[<->,thick,>=stealth] (d) to[bend right=36] (e);
			\draw[<->,thick,>=stealth] (e) to[bend right=36] (a);
			\draw[->,thick,>=stealth] (a) to[bend right=15] (c);
			\draw[->,thick,>=stealth] (b) to[bend right=15] (d);
			\draw[->,thick,>=stealth] (c) to[bend right=15] (e);
			\draw[->,thick,>=stealth] (d) to[bend right=15] (a);
			\draw[->,thick,>=stealth] (e) to[bend right=15] (b);
			\end{tikzpicture}
            \mbox{(3)}
    \end{minipage}&
    \begin{minipage}{0.15\linewidth}
        \centering
        \begin{tikzpicture}[auto,scale=0.7,rotate=0]
			\node[draw,circle,inner sep=2pt] (a) at (0,1){};
			\node[draw,circle,inner sep=2pt] (b) at (sin{-72}, cos{-72}){};
			\node[draw,circle,inner sep=2pt] (c) at (sin{-144}, cos{-144}){};
			\node[draw,circle,inner sep=2pt] (d) at (sin{-216}, cos{-216}){};
			\node[draw,circle,inner sep=2pt] (e) at (sin{-288}, cos{-288}){};
			\draw[<->,thick,>=stealth] (a) to[bend right=36] (b);
			\draw[<->,thick,>=stealth] (b) to[bend right=36] (c);
			\draw[<->,thick,>=stealth] (c) to[bend right=36] (d);
			\draw[<->,thick,>=stealth] (d) to[bend right=36] (e);
			\draw[<->,thick,>=stealth] (e) to[bend right=36] (a);
			\draw[<->,thick,>=stealth] (a) to[bend right=15] (c);
			\draw[<->,thick,>=stealth] (b) to[bend right=15] (d);
			\draw[<->,thick,>=stealth] (c) to[bend right=15] (e);
			\draw[<->,thick,>=stealth] (d) to[bend right=15] (a);
			\draw[<->,thick,>=stealth] (e) to[bend right=15] (b);
			\end{tikzpicture}
            \mbox{(4)}
    \end{minipage}
    \end{tabular}
    \caption{The case $n=10$}
    \label{fig:order10}
\end{figure}
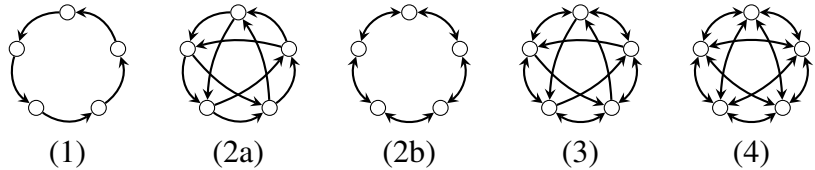
Therefore, we have completed the classification of the case $n=10$.
\end{example}

The next simplest case is when the order can be written as $3p$ by some prime $p$.
\begin{example}[$n=9$]
Since $9$ is the square of $3$, $Q$ has three $\Inn(Q)$-orbits unless $Q$ is trivial.
Thus, possible representatives of lists are $(1,0),(1,1)$, and $(1,2)$ (Figure \ref{fig:order9}).
Moreover, the quandles determined by these weights are not isomorphic to each other.
\begin{figure}[h]
    \centering
    \begin{tabular}{ccc}
    \begin{minipage}{0.2\linewidth}
    \centering
    \begin{tikzpicture}[auto,scale=0.7]
	\node[draw,circle, inner sep=2pt] (a) at (-1,0){};
	\node[draw,circle, inner sep=2pt] (b) at (1,0){};
	\node[draw,circle, inner sep=2pt] (c) at (0,1.73){};

	\draw[->,thick,>=stealth] (a) to[bend right=30] node[auto=right]{\footnotesize $+1$} (b);
	\draw[->,thick,>=stealth] (b) to[bend right=30] (c);
	\draw[->,thick,>=stealth] (c) to[bend right=30] (a);
	\end{tikzpicture}
    \end{minipage}&
    \begin{minipage}{0.2\linewidth}
    \centering
    \begin{tikzpicture}[auto,scale=0.7]
	\node[draw,circle, inner sep=2pt] (a) at (-1,0){};
	\node[draw,circle, inner sep=2pt] (b) at (1,0){};
	\node[draw,circle, inner sep=2pt] (c) at (0,1.73){};

	\draw[->,thick,>=stealth] (a) to[bend right=30] node[auto=right]{\footnotesize $+1$} (b);
	\draw[->,thick,>=stealth] (b) to[bend right=30] (c);
	\draw[->,thick,>=stealth] (c) to[bend right=30] (a);

	\draw[->,thick,>=stealth] (a) to[bend right=15] (c);
	\draw[->,thick,>=stealth] (b) to[bend right=15] (a);
	\draw[->,thick,>=stealth] (c) to[bend right=15] (b);
	\end{tikzpicture}
    \end{minipage}&
    \begin{minipage}{0.2\linewidth}
    \centering
    \begin{tikzpicture}[auto,scale=0.7]
	\node[draw,circle, inner sep=2pt] (a) at (-1,0){};
	\node[draw,circle, inner sep=2pt] (b) at (1,0){};
	\node[draw,circle, inner sep=2pt] (c) at (0,1.73){};

	\draw[->,thick,>=stealth] (a) to[bend right=30] node[auto=right]{\footnotesize $+1$} (b);
	\draw[->,thick,>=stealth] (b) to[bend right=30] (c);
	\draw[->,thick,>=stealth] (c) to[bend right=30] (a);

	\draw[->,thick,dashed,>=stealth] (a) to[bend right=15] node[auto=right]{\footnotesize $-1$} (c);
	\draw[->,thick,dashed,>=stealth] (b) to[bend right=15] (a);
	\draw[->,thick,dashed,>=stealth] (c) to[bend right=15] (b);
	\end{tikzpicture}
    \end{minipage}
    \end{tabular}
    \caption{The case $n=9$}
    \label{fig:order9}
\end{figure}
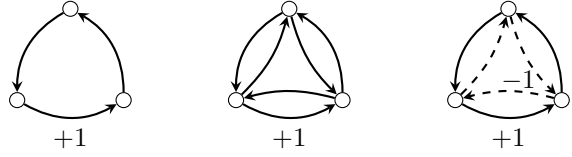
\end{example}

\begin{example}[$n=8$]
The quandles are uniquely determined when $|X|=8$ or $|X|=2$.
When $|X|=4$, then $A=\bZ_2$.
If there exists $f\in H_Q$ such that $\ord f=4$, there are $2^3-1=7$ possible weight lists $(d(x_1,x_0),d(x_2,x_0),d(x_3,x_0))$.
These weights are classified into five types as follows
\begin{itemize}
    \item[(1a)] $(1,0,0),(0,0,1)$;
    \item[(1b)] $(0,1,0)$;
    \item[(2a)] $(1,0,1)$;
    \item[(2b)] $(0,1,1),(1,1,0)$;
    \item[(3)] $(1,1,1)$.
\end{itemize}
Furthermore, they give different support graphs, as shown in Figure \ref{fig:order8}.
    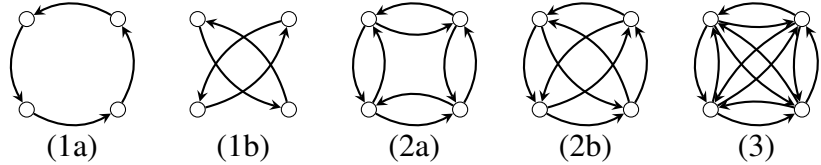
\begin{figure}[h]
    \centering
    \begin{tabular}{ccccc}
    \begin{minipage}{0.15\linewidth}
    \centering
    \begin{tikzpicture}[auto,scale=1.2]
	\node[draw,circle, inner sep=2pt] (a) at (0,0){};
	\node[draw,circle, inner sep=2pt] (b) at (1,0){};
	\node[draw,circle, inner sep=2pt] (c) at (1,1){};
    \node[draw,circle, inner sep=2pt] (d) at (0,1){};
	\draw[->,thick,>=stealth] (a) to[bend right=30] (b);
	\draw[->,thick,>=stealth] (b) to[bend right=30] (c);
	\draw[->,thick,>=stealth] (c) to[bend right=30] (d);
    \draw[->,thick,>=stealth] (d) to[bend right=30] (a);
	\end{tikzpicture}
    \end{minipage}&
    \begin{minipage}{0.15\linewidth}
    \centering
    \begin{tikzpicture}[auto,scale=1.2]
	\node[draw,circle, inner sep=2pt] (a) at (0,0){};
	\node[draw,circle, inner sep=2pt] (b) at (1,0){};
	\node[draw,circle, inner sep=2pt] (c) at (1,1){};
    \node[draw,circle, inner sep=2pt] (d) at (0,1){};
	\draw[->,thick,>=stealth] (a) to[bend right=30] (c);
	\draw[->,thick,>=stealth] (b) to[bend right=30] (d);
	\draw[->,thick,>=stealth] (c) to[bend right=30] (a);
    \draw[->,thick,>=stealth] (d) to[bend right=30] (b);
	\end{tikzpicture}
    \end{minipage}&
    \begin{minipage}{0.15\linewidth}
    \centering
    \begin{tikzpicture}[auto,scale=1.2]
	\node[draw,circle, inner sep=2pt] (a) at (0,0){};
	\node[draw,circle, inner sep=2pt] (b) at (1,0){};
	\node[draw,circle, inner sep=2pt] (c) at (1,1){};
    \node[draw,circle, inner sep=2pt] (d) at (0,1){};
	\draw[->,thick,>=stealth] (a) to[bend right=30] (b);
	\draw[->,thick,>=stealth] (b) to[bend right=30] (c);
	\draw[->,thick,>=stealth] (c) to[bend right=30] (d);
    \draw[->,thick,>=stealth] (d) to[bend right=30] (a);
    \draw[->,thick,>=stealth] (a) to[bend right=30] (d);
	\draw[->,thick,>=stealth] (d) to[bend right=30] (c);
	\draw[->,thick,>=stealth] (c) to[bend right=30] (b);
    \draw[->,thick,>=stealth] (b) to[bend right=30] (a);
	\end{tikzpicture}
    \end{minipage}&
    \begin{minipage}{0.15\linewidth}
    \centering
    \begin{tikzpicture}[auto,scale=1.2]
	\node[draw,circle, inner sep=2pt] (a) at (0,0){};
	\node[draw,circle, inner sep=2pt] (b) at (1,0){};
	\node[draw,circle, inner sep=2pt] (c) at (1,1){};
    \node[draw,circle, inner sep=2pt] (d) at (0,1){};
	\draw[->,thick,>=stealth] (a) to[bend right=30] (b);
	\draw[->,thick,>=stealth] (b) to[bend right=30] (c);
	\draw[->,thick,>=stealth] (c) to[bend right=30] (d);
    \draw[->,thick,>=stealth] (d) to[bend right=30] (a);
    \draw[->,thick,>=stealth] (a) to[bend right=30] (c);
	\draw[->,thick,>=stealth] (b) to[bend right=30] (d);
	\draw[->,thick,>=stealth] (c) to[bend right=30] (a);
    \draw[->,thick,>=stealth] (d) to[bend right=30] (b);
	\end{tikzpicture}
    \end{minipage}&
    \begin{minipage}{0.15\linewidth}
    \centering
    \begin{tikzpicture}[auto,scale=1.2]
	\node[draw,circle, inner sep=2pt] (a) at (0,0){};
	\node[draw,circle, inner sep=2pt] (b) at (1,0){};
	\node[draw,circle, inner sep=2pt] (c) at (1,1){};
    \node[draw,circle, inner sep=2pt] (d) at (0,1){};
	\draw[->,thick,>=stealth] (a) to[bend right=30] (b);
	\draw[->,thick,>=stealth] (b) to[bend right=30] (c);
	\draw[->,thick,>=stealth] (c) to[bend right=30] (d);
    \draw[->,thick,>=stealth] (d) to[bend right=30] (a);
    \draw[->,thick,>=stealth] (a) to[bend right=15] (d);
	\draw[->,thick,>=stealth] (d) to[bend right=15] (c);
	\draw[->,thick,>=stealth] (c) to[bend right=15] (b);
    \draw[->,thick,>=stealth] (b) to[bend right=15] (a);
    \draw[->,thick,>=stealth] (a) to[bend right=15] (c);
	\draw[->,thick,>=stealth] (b) to[bend right=15] (d);
	\draw[->,thick,>=stealth] (c) to[bend right=15] (a);
    \draw[->,thick,>=stealth] (d) to[bend right=15] (b);
	\end{tikzpicture}
    \end{minipage}\\
    \mbox{(1a)}&\mbox{(1b)}&\mbox{(2a)}&\mbox{(2b)}&\mbox{(3)}
    \end{tabular}
    \caption{The case $n=8$}
    \label{fig:order8}
\end{figure}
If there is no $f\in H_Q$ such that $\ord f=4$, $H_Q$ contains permutations $(x_0x_1)(x_2x_3)$ and $(x_0x_2)(x_1x_3)$ on $X=\{x_0,x_1,x_2,x_3\}$ up to isomorphisms.
However, graphs that are equivalent under the action of $H_Q$ are just the above graphs.
\end{example}

Even in the case of larger orders, if the order is written as the product of two primes $p$ and $q$, then we can apply Burnside's lemma. 
It is because the isomorphism classes are determined by $\Aut(\bZ_p)=\bZ_p^\times$ and $\Aut(\bZ_q)=\bZ_q^\times$ actions.

Let us recall Burnside's lemma.
\begin{lemma}[Burnside's lemma \cite{ste}]
    Let $G$ be a group and let $\ord\sigma$ denote the order of $\sigma$ for $\sigma\in G$.
    Suppose $G$ acts on a finite set $\Omega$.
    The number of $G$-orbit of $\Omega$ is given by
    \begin{equation*}
        \frac{1}{|G|}\sum_{\sigma\in G}|\Omega^\sigma|
    \end{equation*}
    where $\Omega^\sigma=\{\omega\in\Omega\mid \sigma(\omega)=\omega\}$.
\end{lemma}
Let us calculate some examples using this lemma.
\begin{example}[$n=10,14,15$]\mbox{}
\begin{itemize}
\item[$n=10$:] Consider the set 
    \begin{equation*}
        \Omega_{\bZ_5,\bZ_2}=\{(d(x_i,x_0))_{i\in \bZ_5^\times}\mid d(x_i,x_0)\in \bZ_2(i=1,\dots,4)\}
    \end{equation*}
of all lists of weight is acted by $\bZ_5^\times\cong \bZ_4$.
We get $|\Omega_{\bZ_5,\bZ_2}|=2^{(5-1)}$.
If $\sigma\in\mathbb{Z}^\times_5$ is of order $k$, there are $(5-1)/k$ orbits of the action $\sigma$ on $\Omega_{\bZ_5,\bZ_2}$.
Notice that the all zero list $(0,\dots,0)$ gives the trivial quandle.
Adding one, which is the case $|X|=2$, we have that the number of equivalent classes of homogeneous quandles of order $10$ with abelian inner automorphisms is equal to
\begin{align*}
    &1+1+(|\Omega_{\bZ_5,\bZ_2}/\bZ_5^\times|-1)\\
    =&1+\frac{1}{5-1}\sum_{\sigma\in \bZ_5^\times}|\Omega_{\bZ_5,\bZ_2}^\sigma|\\
    =&1+\frac{1}{4}\left(2^{4/\ord(1)}+2^{4/\ord(2)}+2^{4/\ord(3)}+2^{4/\ord(4)}\right)\\
    =&1+\frac{1}{4}\left(2^{4/1}+2^{4/4}+2^{4/4}+2^{4/2}\right)\\
    =&7=1+1+5.
\end{align*}
This result coincides with the enumeration in Example \ref{exam:n=10}.

\item[$n=14$:] Similarly, let $\Omega_{\bZ_7,\bZ_2}$ denote
\begin{equation*}
        \{(d(x_i,x_0))_{i\in \bZ_7^\times}\mid d(x_i,x_0)\in \bZ_2(i=1,\dots,6)\}.
\end{equation*}
We have that the number of equivalent classes of homogeneous quandles of order $14$ with abelian inner automorphisms is equal to
\begin{align*}
    &1+1+(|\Omega_{\bZ_7,\bZ_2}/\bZ_7^\times|-1)\\
    =&1+\frac{1}{7-1}\sum_{\sigma\in \bZ_7^\times}|\Omega_{\bZ_7,\bZ_2}^\sigma|\\
    =&1+\frac{1}{6}\left(2^\frac{6}{\ord(1)}+2^\frac{6}{\ord(2)}+2^\frac{6}{\ord(3)}+2^\frac{6}{\ord(4)}+2^\frac{6}{\ord(5)}+2^\frac{6}{\ord(6)}\right)\\
    =&1+\frac{1}{6}\left(2^{6/1}+2^{6/3}+2^{6/6}+2^{6/3}+2^{6/6}+2^{6/2}\right)\\
    =&15.
\end{align*}
\item[$n=15$:] In this case, we need to count when $|X|=3$ or when $|X|=5$ separately.
Let $\Omega_{\bZ_p,\bZ_q}$ denote
\begin{equation*}
        \{(d(x_i,x_0))_{i\in \bZ_p^\times}\mid d(x_i,x_0)\in \bZ_q(i=1,\dots,p-1)\}.
\end{equation*}
Then, the number of equivalent classes of homogeneous quandles of order $15$ with abelian inner automorphisms is equal to
\begin{align*}
    &1+(|\Omega_{\bZ_3,\bZ_5}/\bZ_3^\times\times \bZ_5^\times|-1)+(|\Omega_{\bZ_5,\bZ_3}/\bZ_5^\times\times\bZ_3^\times|-1)\\
    =&1+\left(\frac{1}{(3-1)(5-1)}\sum_{\sigma\in \bZ_3^\times\times \bZ_5^\times}|\Omega_{\bZ_3,\bZ_5}^\sigma|-1\right)\\&+\left(\frac{1}{(5-1)(3-1)}\sum_{\sigma\in \bZ_5^\times\times \bZ_3^\times}|\Omega_{\bZ_5,\bZ_3}^\sigma|-1\right)\\
    =&1+(40/8-1)+(112/8-1)=1+4+13=18.
\end{align*}
\end{itemize}
\end{example}

As in the case of order $10=2\cdot5,14=2\cdot7$, we obtain the following counting formula for quandles of order $2p$.
\begin{proposition}
    Let $p$ be a prime number.
    Then, the number of isomorphism classes of the homogeneous quandles with abelian inner automorphism groups of order $2p$ is given by
    \begin{equation*}
        1+\frac{1}{p-1}\sum_{i=1}^{p-1}2^{(p-1)/\gcd(p-1,i)}.
    \end{equation*}
\end{proposition}
Here we used $\ord i=\gcd(p-1,i)$ for $i\in (\bZ_{p-1},+)$.
This formula suggests that the number of isomorphism classes of $X\times_d A$ of order $2p$ increases by about $2^p/p$.

\section{Remarks on Theorem \ref{thm:main_transitive}}\label{sec:rema}
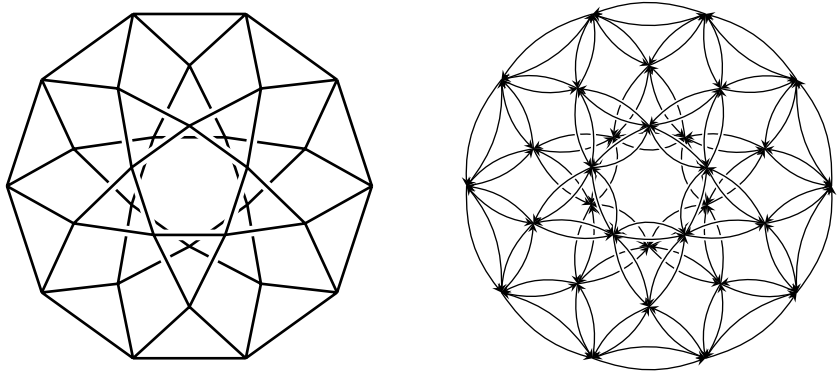
\begin{figure}[h]
\centering
\begin{tabular}{cc}
\begin{minipage}{0.45\linewidth}
\centering
\begin{tikzpicture}[scale=0.8]
\foreach \arg in {0,72,144,216,288}{
	\draw[line width = 1pt, rotate around={\arg+36:(0,0)},color=black, line join=bevel] (3,0)--(2*cos{-18}, 2*sin{-18})--(cos{18}, sin{18})--( cos{-54}, sin{-54})--(2*cos{-18}, 2*sin{-18})--(3*cos{-36}, 3*sin{-36})--(3,0);}
\foreach \arg in {0,72,144,216,288}{
	\draw[line width = 4pt, rotate around={\arg:(0,0)},color=white] (2*cos{-18}, 2*sin{-18})--( cos{18}, sin{18})--( cos{-54}, sin{-54})--(2*cos{-18}, 2*sin{-18});}
\foreach \arg in {0,72,144,216,288}{
	\draw[line width = 1pt, rotate around={\arg:(0,0)},color=black, line join=bevel] (3,0)--(2*cos{-18}, 2*sin{-18})--(cos{18}, sin{18})--( cos{-54}, sin{-54})--(2*cos{-18}, 2*sin{-18})--(3*cos{-36}, 3*sin{-36})--(3,0);}
\end{tikzpicture}
\end{minipage}&
\begin{minipage}{0.45\linewidth}
\centering
\begin{tikzpicture}[scale=0.8]
\foreach \arg in {0,72,144,216,288}{
	\draw[line width = .5pt, rotate around={\arg+36:(0,0)},color=black, arrows=-stealth] (3,0)  to [bend right=20](2*cos{-18}, 2*sin{-18});
	\draw[line width = .5pt, rotate around={\arg+36:(0,0)},color=black, arrows=stealth-] (3,0) to [bend right=-20](2*cos{-18}, 2*sin{-18});
	\draw[line width = .5pt, rotate around={\arg+36:(0,0)},color=black, arrows=-stealth] (2*cos{-18}, 2*sin{-18}) to [bend right=20]( cos{18}, sin{18});
	\draw[line width = .5pt, rotate around={\arg+36:(0,0)},color=black, arrows=stealth-] (2*cos{-18}, 2*sin{-18}) to [bend right=-20](cos{18}, sin{18});
	\draw[line width = .5pt, rotate around={\arg+36:(0,0)},color=black, arrows=-stealth] (cos{18}, sin{18}) to [bend right=36](cos{-54}, sin{-54});
	\draw[line width = .5pt, rotate around={\arg+36:(0,0)},color=black, arrows=stealth-] (cos{18}, sin{18}) to [bend right=-12](cos{-54}, sin{-54});
	\draw[line width = .5pt, rotate around={\arg+36:(0,0)},color=black, arrows=-stealth] (cos{-54}, sin{-54}) to [bend right=20](2*cos{-18}, 2*sin{-18});
	\draw[line width = .5pt, rotate around={\arg+36:(0,0)},color=black, arrows=stealth-] (cos{-54}, sin{-54}) to [bend right=-20](2*cos{-18}, 2*sin{-18});
	\draw[line width = .5pt, rotate around={\arg+36:(0,0)},color=black, arrows=-stealth] (2*cos{-18}, 2*sin{-18}) to [bend right=20](3*cos{-36}, 3*sin{-36});
	\draw[line width = .5pt, rotate around={\arg+36:(0,0)},color=black, arrows=stealth-] (2*cos{-18}, 2*sin{-18}) to [bend right=-20](3*cos{-36}, 3*sin{-36});
	\draw[line width = .5pt, rotate around={\arg+36:(0,0)},color=black, arrows=-stealth] (3,0) to [bend right=20](3*cos{-36}, 3*sin{-36});
	\draw[line width = .5pt, rotate around={\arg+36:(0,0)},color=black, arrows=stealth-] (3,0) to [bend right=-20](3*cos{-36}, 3*sin{-36});}
\foreach \arg in {0,72,144,216,288}{
	\draw[line width = 2pt, rotate around={\arg:(0,0)},color=white] (2*cos{-18}, 2*sin{-18}) to [bend right=20](cos{18}, sin{18});
	\draw[line width = 2pt, rotate around={\arg:(0,0)},color=white] (2*cos{-18}, 2*sin{-18}) to [bend right=-20](cos{18}, sin{18});
	\draw[line width = 2pt, rotate around={\arg:(0,0)},color=white] (cos{18}, sin{18}) to [bend right=36](cos{-54}, sin{-54});
	\draw[line width = 2pt, rotate around={\arg:(0,0)},color=white] (cos{18}, sin{18}) to [bend right=-12](cos{-54}, sin{-54});
	\draw[line width = 2pt, rotate around={\arg:(0,0)},color=white] (cos{-54}, sin{-54}) to [bend right=20](2*cos{-18}, 2*sin{-18});
	\draw[line width = 2pt, rotate around={\arg:(0,0)},color=white] (cos{-54}, sin{-54}) to [bend right=-20](2*cos{-18}, 2*sin{-18});}
\foreach \arg in {0,72,144,216,288}{
	\draw[line width = .5pt, rotate around={\arg:(0,0)},color=black, arrows=-stealth] (3,0) to [bend right=20](2*cos{-18}, 2*sin{-18});
	\draw[line width = .5pt, rotate around={\arg:(0,0)},color=black, arrows=stealth-] (3,0) to [bend right=-20](2*cos{-18}, 2*sin{-18});
	\draw[line width = .5pt, rotate around={\arg:(0,0)},color=black, arrows=-stealth] (2*cos{-18}, 2*sin{-18}) to [bend right=20]( cos{18}, sin{18});
	\draw[line width = .5pt, rotate around={\arg:(0,0)},color=black, arrows=stealth-] (2*cos{-18}, 2*sin{-18}) to [bend right=-20](cos{18}, sin{18});
	\draw[line width = .5pt, rotate around={\arg:(0,0)},color=black, arrows=-stealth] (cos{18}, sin{18}) to [bend right=36](cos{-54}, sin{-54});
	\draw[line width = .5pt, rotate around={\arg:(0,0)},color=black, arrows=stealth-] (cos{18}, sin{18}) to [bend right=-12](cos{-54}, sin{-54});
	\draw[line width = .5pt, rotate around={\arg:(0,0)},color=black, arrows=-stealth] (cos{-54}, sin{-54}) to [bend right=20](2*cos{-18}, 2*sin{-18});
	\draw[line width = .5pt, rotate around={\arg:(0,0)},color=black, arrows=stealth-] (cos{-54}, sin{-54}) to [bend right=-20](2*cos{-18}, 2*sin{-18});
	\draw[line width = .5pt, rotate around={\arg:(0,0)},color=black, arrows=-stealth] (2*cos{-18}, 2*sin{-18}) to [bend right=20](3*cos{-36}, 3*sin{-36});
	\draw[line width = .5pt, rotate around={\arg:(0,0)},color=black, arrows=stealth-] (2*cos{-18}, 2*sin{-18}) to [bend right=-20](3*cos{-36}, 3*sin{-36});
	\draw[line width = .5pt, rotate around={\arg:(0,0)},color=black, arrows=-stealth] (3,0) to [bend right=20](3*cos{-36}, 3*sin{-36});
	\draw[line width = .5pt, rotate around={\arg:(0,0)},color=black, arrows=stealth-] (3,0) to [bend right=-20](3*cos{-36}, 3*sin{-36});}
\end{tikzpicture}
\end{minipage}
\end{tabular}
\caption{Left: the $1$-skeleton $\tilde{\Gamma}_{ID}$ of the icosidodecahedron, right: the symmetric digraph of $\tilde{\Gamma}_{ID}$}\label{fig:IDgraph}
\end{figure}

The converse of Theorem \ref{thm:main_transitive} does not hold in general.
This means that even if $X$ is homogeneous and what isomorphism $\Inn(Q)/\Inn(Q)_x \cong \{x\}\times A$ we choose for each $x\in X$, $(X,d)$ may not become $(X,d)$-homogeneous.
We provide such an example as follows(see Figure \ref{fig:IDgraph}).

\begin{figure}[h]
    \centering
    \begin{tikzpicture}[xscale=2,auto,->]
        \node (x1) at (-1,0.5) {$x_1$};\node (x2) at (1,0.5) {$x_2$};\node (x3) at (1,-0.5) {$x_3$};\node (x4) at (-1,-0.5) {$x_4$};\node (y) at (0,0) {$y$};
        \draw (x1) to node {\scriptsize{$+2^0$}} (y);\draw (x2) to node[swap] {\scriptsize{$+2^1$}} (y);\draw (x3) to node {\scriptsize{$+2^2$}} (y);\draw (x4) to node[swap] {\scriptsize{$+2^3$}} (y);
    \end{tikzpicture}
    \caption{$d(x_i,y)$ of $Q_{ID}$}\label{fig:weightID}
\end{figure}
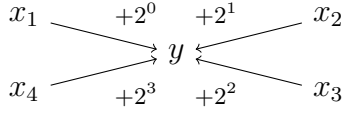

Let $X$ be the vertex set of the icosidodecahedron, and let $A = \bZ_5$.
For each vertex $y$, there are just four vertices adjacent to $y$. We put them by $x_1, x_2, x_3$, and $x_4$ in the clockwise.
Take the non-zero elements of $A$ for ordered pairs $(x_i,y)$ $(i=1,2,3,4)$ such that they satisfy the following condition (Figure \ref{fig:weightID});
\begin{equation}\label{eq:condID}
    d(x_i,y)=2^{i-1}d(x_1,y).
\end{equation}

The $\bZ_5$-weighted graph $(X,d)$ is indecomposable since 
\begin{equation}\label{eq:arrows}
\{d(x,y)\in\bZ_5\mid x\in X\}=\bZ_5    
\end{equation}
for each $y\in X$.
Remark that the automorphism group of $\bZ_5$ is generated by the map $\times 2$.
If $(X,d')$ satisfies the same condition \eqref{eq:condID}, taking $\sigma:X\to \Aut(\bZ_5)$ by $y\mapsto \sigma_y=\times d(x,y)/d'(x,y)\in \Aut(\bZ_5)$, then $(\mathrm{id},\sigma):X\times_d \bZ_5\to X\times_{d'} \bZ_5$ becomes an isomorphism by Proposition \ref{prop:flip}.
Thus, the quandle $X\times_{d} A$ does not depend on the choice of $d(x_1,y)$ for each $y\in X$ up to isomorphism.
\begin{definition}
    Define the quandle $Q_{ID}$ by $X\times_{d} A$ in the above setting.
\end{definition}

    Let $\tilde{\Gamma}_{ID}$ be the $1$-skeleton of the icosidodecahedron and consider $\tilde{\Gamma}$ as an undirected graph. Denote by $\Gamma_{ID}$ the symmetric directed graph obtained by replacing the edges $\{x,y\}$ of $\tilde{\Gamma}_{ID}$ with two directed edges $(x,y)$ and $(y,x)$.

\begin{proposition}
    The quandle $Q_{ID}$ is homogeneous.
\end{proposition}
\proof
By the construction, we obtain the graph isomorphism $\Gamma(X,d)\cong \Gamma_{ID}$.
The graph $\tilde{\Gamma}_{ID}$ is a vertex-transitive graph by the symmetry group $\Aut(\Gamma_{ID})$ of the icosidodecahedron.
In particular, the subgroup $\Aut(\Gamma_{ID})^{+}$ consisting of  transformations preserving the orientation of the icosidodecahedron also acts transitively on $\Gamma_{ID}$.

A symmetric transformation $g\in \Aut(\Gamma_{ID})^{+}$ that preserves the orientation may give a different $\bZ_5$-weight $d'=d\circ (g,g)$.
That is, $d(x,y)= d'(x,y)=d(g(x),g(y))$ does not hold in general.
However, this $d'$ satisfies condition \eqref{eq:condID} again.
Therefore, Proposition \ref{prop:flip} yields that the quandle $Q_{ID}$ is homogeneous.
\endproof

Note that the automorphism of the $1$-skeleton of the polytope induces the automorphism of the face lattice.
Naturally, we consider that $\Aut(X,d)$ acts on $\Gamma_{ID}$ and the edges and faces of the icosidodecahedron.

\begin{figure}[h]
    \centering
    \begin{tabular}{ccc}
    \begin{minipage}{0.3\linewidth}
    \centering
    \begin{tikzpicture}[scale=1.5]
    \foreach \arg in {0,72,144,216,288}{\node[inner sep=0.2em, fill=black!100, circle] (\arg) at (cos{(\arg+90)}, sin{(\arg+90)}) {};\node[coordinate] (a\arg) at (cos{(\arg+18)}, sin{(\arg+18)}) {};\node[coordinate] (b\arg) at (cos{(\arg+162)}, sin{(\arg+162)}) {};}
    \foreach \arg in {0,72,144,216,288}{
    \draw[->,thick] (a\arg) to[bend right=18] (\arg);\draw[->,dashed,thick] (b\arg) to[bend right=18] (\arg);}
    \node (A) at (0,0){$F_1$};
    \node (d+) at (0,-0.5){$d^+$};
    \node (d-) at (0,-1.1){$d^-$};
    \end{tikzpicture}
    \end{minipage}&
    \begin{minipage}{0.3\linewidth}
    \centering
    \begin{tikzpicture}[scale=0.7]
    \foreach \arg in {0,120,240}{
	\draw[line width = 1pt, rotate around={\arg:(0,0)},line join=bevel] (cos{-30}, sin{-30})--(cos{90}, sin{90})--(2.5*cos{60}, 2.5*sin{60})--(3*cos{30}, 3*sin{30})--(2.5,0)--cycle;}
    \node (A) at (1.5*cos{30}, 1.5*sin{30}){$F_1$};\node (B) at (1.5*cos{150}, 1.5*sin{150}){$F_3$};\node (C) at (1.5*cos{270}, 1.5*sin{270}){$F_2$};\node (D) at (0,0){$F_0$};\node (y) at (1.5*cos{210}, 1.5*sin{210}){$y$};\node (x0) at (1.5*cos{330}, 1.5*sin{330}){$x_0$};\node (x1) at (1.5*cos{90}, 1.5*sin{90}){$x_1$};
    \end{tikzpicture}
    \end{minipage}&
    \begin{minipage}{0.3\linewidth}
    \centering
    \begin{tikzpicture}[yscale=1.5]
    \node[inner sep=0.2em, fill=black!100, circle] (a0) at (0,0) {};
    \foreach \int in {1,2,3,4}{\node[inner sep=0.2em, fill=black!100, circle] (a\int) at (cos{(\int*90-45)},sin{(\int*90-45)}) {};}
    \foreach \int in {1,3}{\draw[->,thick] (a\int) to[bend right=18] (a0);\draw[<-,dashed,thick] (a\int) to[bend right=-18] (a0);}
    \foreach \int in {2,4}{\draw[->,dashed,thick] (a\int) to[bend right=18] (a0);\draw[<-,thick] (a\int) to[bend right=-18] (a0);}
    \node (z) at (0,0.3) {$z$};\node (F1) at (-1,0) {$F_1$};\node (gF1) at (1,0) {$g(F_1)$};
    \end{tikzpicture}
    \end{minipage}
    \end{tabular}
    \caption{Vertices, weighted edges, and faces in Proposition \ref{prop:ID}}\label{fig:PfOfID}
\end{figure}
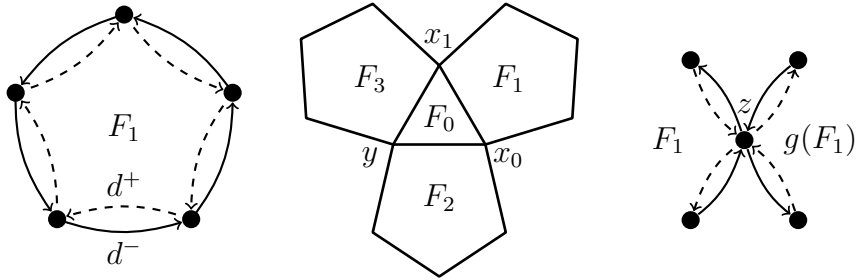

\begin{proposition}\label{prop:ID}
    There is no homogeneous $\bZ_5$-weight $d:X\times X\to \bZ_5$ on $X$ such that $X\times_d A$ is isomorphic to $Q_{ID}$.
\end{proposition}
\proof
Suppose $\Aut(X,d)$ acts transitively on $X$.
By assumption and $|X|=2\times 3\times 5$, there exists an automorphism $f:X\to X$ of order $5$.
Since $f$ induces the automorphism of the support graph $\Gamma(X,d)$, this is a rotation of the icosidodecahedron, which fixes some pentagon face $F_1$.
We denote the vertices of $F_1$ by $x_0,\dots,x_4$ clockwise.
Then there exist $d^-,d^+\in \bZ_5\setminus \{0\}$ such that $d(v_i,v_{i-1})=d^-$ and $d(v_i,v_{i+1})=d^+$ for each $i$ because $f$ preserves $A$-weight (see the left in Figure \ref{fig:PfOfID}).
Let $F_0$ be the triangle that has $(x_0,x_1)$ as an edge, and denote by $y$ the other vertex of $F_0$ (see the center in Figure \ref{fig:PfOfID}).
By homogeneity of $(X,d)$, there exists $g\in \Aut(X,d)$ such that $g(x_0)=y$.
Let $F_2$ and $F_3$ be the other pentagons having a common edge with $F_0$.
Then the image $g(F_1)$ of $F_1$ through $g$ is $F_2$ or $F_3$.
Then $g(F_1)$ has a common vertex $z\in\{x_0,x_1\}$ with $F_1$ (see the right in Figure \ref{fig:PfOfID}).
Because $g$ preserves $A$-weight, we have $|\{d(x,z)\in \bZ_5\mid x\in X\}|=3$.
This contradicts the equation \eqref{eq:arrows}.

\endproof


\begin{thebibliography}{99}

\bibitem{bar-nas}
V. Bardakov and T. Nasybullov, Embeddings of quandles into groups, \emph{J. Algebra Appl.} {\bf 19} (2020), no. 7, 2050136, 20 pp.

\bibitem{cens}
J. S. Carter, M. Elhamdadi, M. A. Nikiforou, M. Saito, Extensions of quandles and cocycle knot invariants, \emph{J. Knot Theory Ramifications}  {\bf 12} (2003), no. 6, 725–738.


\bibitem{cks}
J. S.Carter, S. Kamada, M. Saito, Diagrammatic computations for quandles and cocycle, \emph{Diagrammatic Morphisms and Applications San Francisco, CA, 2000}, 51-74. Contemp. Math., 318 Amer. Math. Soc., Providence, RI, 2003.

\bibitem{fur-tam}
K. Furuki, H. Tamaru, Homogeneous quandles with abelian inner automorphism groups and vertex-transitive graphs, preprint, arXiv:2403.06209.


\bibitem{hig-kur}
A. Higashitani, H. Kurihara, Homogeneous quandles arising from automorphisms of symmetric groups,
\emph{Comm. Algebra} {\bf 51} (2023), no.4, 1413–1430.

\bibitem{ish-tam}
Y. Ishihara, H. Tamaru, Flat connected finite quandles, \emph{Proc. Amer. Math. Soc.}, {\bf 144} (2016), no. 11, 4959-4971.

\bibitem{joy}
D. Joyce, A classifying invariant of knots, the knot quandle, 
\emph{J. Pure Appl. Algebra}, {\bf 23} (1982), no. 1, 37-65.


\bibitem{jpsz}
P. Jedli\v{c}ka, A. Pilitowska, D. Stanovsk\'{y}, A. Zamojska-Dzienio, The structure of medial quandles, \emph{J. Algebra}, {\bf 443} (2015), 300–334.


\bibitem{kmtt}
A. Kubo, M. Nagashiki, T. Okuda, H. Tamaru, A commutativity condition for subsets in quandles—a generalization of antipodal subsets,
\emph{Differential Geometry and Global Analysis, in Honor of Tadashi Nagano}, 103-125. Contemp. Math., 777, Amer. Math. Soc., Providence RI, 2022.


\bibitem{leb-mor}
V. Lebed and A. Mortier, Abelian quandles and quandles with abelian structure group, 
\emph{J. Pure Appl. Algebra}, {\bf 225}, (2021),no. 1, Paper No. 106474, 22 pp.

\bibitem{mat}
S. V. Matveev, Distributive groupoids in knot theory, \emph{Mat. Sb. (N.S.)} {\bf 119} (1982), no. 1, 78–88.

\bibitem{muz}
M. Muzychuk. On \'Ad\'am’s conjecture for circulant graphs. 
\emph{Discrete Math.} {\bf 176} (1997), no. 1-3, 285–298.

\bibitem{sin}
M. Singh, Classification of flat connected quandles, \emph{J. Knot Theory Ramifications} {\bf 25} (2016), no. 13, 1650071, 8 pp.

\bibitem{ste}
B. Steinberg, Representation theory of finite groups. An introductory approach. Universitext. \emph{Springer, New York}, (2012) xiv+157 pp.

\bibitem{tak}
M. Takasaki, Abstraction of symmetric transformations. 
\emph{T\^{o}hoku Math. J.} {\bf 49} (1943), 145–207. 


\bibitem{tam}
H. Tamaru, Two-point homogeneous quandles with prime cardinality,
\emph{J. Math. Soc. Japan}, {\bf 65} (2013), no. 4, 1117-1134.
\end{thebibliography}
\end{document}